\documentclass{article}
\usepackage{latexsym}
\usepackage{color}
\usepackage{amssymb}
\usepackage{amsthm}

\usepackage{amssymb}

\usepackage{mathrsfs}

\usepackage{anysize}
\marginsize{2.5cm}{2.5cm}{2.5cm}{3.4cm}

\newtheorem{defi}{Definition}[section]
\newtheorem{prop}[defi]{Proposition}
\newtheorem{teo}[defi]{Theorem}

\newtheorem{cor}[defi]{Corollary}

\def\Qco{\mathfrak{Qcoh}}

\def\Qcart{QMod_{\rm cart}(R)}
\def\cart{Mod_{\rm cart}(R)}

\def\O{{\mathcal O}}

\def\C{{\mathscr C}}

\def\Ker{\mbox{\rm Ker }}

\def\Z{\mathbb{Z}}

\def\fiz{\leftarrow}

\begin{document}

\title{Cartesian modules in small categories}

\bigskip
\author{ \\
Edgar Enochs \\ enochs@ms.uky.edu \\
Department of Mathematics \\
University of Kentucky \\ Lexington, KY 40506-0027 \\ U.S.A \\ \\
\\
\ \ \ Sergio Estrada  \\ \ \ sestrada@um.es \\
\ \
Departamento de Matem\'atica Aplicada \\
\ \ \ Universidad de Murcia \\
\ \ 30100 Murcia \\
SPAIN}

\date{}

\maketitle

\bigskip

\thispagestyle{empty}

\begin{abstract}
In this note we extend the main results of \cite{EE} to the category of
cartesian modules over a flat presheaf of rings $R$ and on an arbitrary small
category. This provides with new applications of that paper to the categories of
quasi--coherent sheaves on an Artin stack or on a Deligne-Mumford stack.

\end{abstract}



\baselineskip=22pt
\section{Introduction}

In \cite{EE} we develop a method for finding a family of generators of the
so-called category of quasi-coherent $R$-modules 
on an arbitrary quiver (cf. \cite[Corollary 3.5]{EE}) and we prove that the
class of flat quasi--coherent
$R$-modules is covering (cf. \cite[Theorem 4.1]{EE}). The present note is
devoted to showing that the same
arguments of \cite{EE} can also be used in a much more general setup, that is
that of cartesian $R$-modules on a flat presheaf of rings $R$ over a small
category $\C$. This extends the main application in \cite{EE} to the category
$\Qco(X)$ of quasi--coherent sheaves on a scheme $X$ and also to the
category $\Qco({\mathcal X})$ of quasi--coherent sheaves on an Artin stack or
on a Deligne-Mumford stack.

 The authors would like to thank Jon Pridham for pointing out the
possibility of these new applications. We also thank Martin Brandenburg
 for pointing out a possible misunderstanding in the wording of \cite[Section
2]{EE}. 

\section{Cartesian modules on quivers}
A quiver $Q$ is a
directed graph. An edge of
a quiver from a vertex $v_1$ to a vertex $v_2$ is denoted by
$a:v_1\to v_2$ or $v_1\stackrel{a}{\to}v_2$, the symbol $E$ will
denote the set of edges. A quiver $Q$ may be thought as a category
in which the objects are the vertices of $Q$ and the morphisms are
the paths (a path is a sequence of edges) of $Q$. The set of all
vertices will be denoted by $V$.

Let $Q=(V,E)$ be a quiver and let $R$ be a presheaf from $Q$
in the category of commutative rings, that is, for each vertex $v\in V$ we
have a ring $R(v)$ and for an edge $a:v\to w$ we have a ring
homomorphism $R(a^{op}):R(w)\to R(v)$.

We shall say that we have an $R$-module $M$ when we have an
$R(v)$-module $M(v)$ and a morphism $M(a^{op}):M(w)\to M(v)$ for each
edge $a:v\to w$ that is $R(v)$-linear. The $R$-module $M$ is said to be a
cartesian $Q$-module if for each edge $a:v\to w$ as above the morphism
$$R(v)\otimes_{R(w)}M(w)\to M(v)$$ given by $r_v\otimes m_w\mapsto
r_v M(a^{op})(m_w)$, $r_v\in R(v),\ m_w\in M(w)$ is an
$R(v)$-isomorphism.

The category of cartesian $Q$-modules is abelian when $R$ is
such that for an edge $v\to w$, $R(v)$ is a {\sl flat}
$R(w)$-module (so the kernel of a morphism between two cartesian $Q$-modules is
also cartesian). In this case
we say $R$ is flat. Coproducts and colimits may be computed
componentwise so direct limits are exact and, as a result of
Proposition \ref{nocont}, we can find a system of generators in
the category. Therefore, the category of cartesian $Q$-modules is indeed a
Grothendieck category when $R$ is flat.

By the tensor product, $M\otimes_R N$, where $M$ is a right
$R$-module and $N$ a left $R$-module, we mean the $\Z$-module
($\Z(v)=\Z$, for all $v\in V$ and $\Z(a)=id_{\Z}$ for all $a\in
E$) such that
$$(M\otimes_R N)(v)=M(v)\otimes_{R(v)} N(v),$$with $(M\otimes_R N)(a)$
the obvious map. We then get the notion of a flat $R$-module and a flat
cartesian $Q$-module.

Given an arbitrary quiver $Q$ and a flat presheaf of rings $R$ over $Q$, we will
denote by $QMod_{\rm cart}(R)$ the category of cartesian $Q$-modules over $R$. 
\section{Cartesian modules on small categories}

Now let $\C$ be any small category, and let $R$ be a flat presheaf of rings on
$\C$. We will consider the category $\cart$ of cartesian $R$-modules. This is
an abelian category and, as a consequence of Proposition \ref{nocont}, it will
be a Grothendieck category. There is a notion of flat cartesian module as
before. Let $Q$ be the quiver whose vertices are the objects of $\C$, and whose
edges are the morphisms of $\C$. It is then clear that the
category $\cart$ is a full subcategory of the category $\Qcart$. Furthermore it
is also clear that if $M\subseteq N$ in $\Qcart$ and $N$ is a cartesian
$R$-module, then $M$ will be automatically a cartesian $R$-module as well. This
easy observation is
crucial in proving our main result and giving our main applications.

\section{Generators and flat covers in $\cart$}

With the observations made in the previous sections, we can use both Proposition
3.3 and Theorem 4.1 of \cite{EE} to infer that $\cart$ is a Grothendieck
category
admitting flat
covers. Throughout this section we will assume that $\C$ is a small category
and $R$ is a flat presheaf of rings on it. We shall denote by $Q$ the
quiver associated to $\C$.

\begin{defi} Let $M$ be a cartesian $Q$-module. The
cardinality of $M$ is defined as the cardinality of the coproduct
(in the category of sets) of all modules associated to the
vertices $v\in V$, that is
$$|M|=|\sqcup_{v\in V}M(v)|$$
\end{defi}

\begin{prop}\label{nocont}

Let $\C$ be any small category with associated quiver $Q_{\C}=(V,E)$
and $M$ a cartesian $R$-module. Let
$\kappa$ be an infinite cardinal such that $\kappa\geq |R(v)|$ for
all $v$ and such that $\kappa\geq max\{|E|,|V|\}$. Let
$X_v\subseteq M(v)$ be subsets with $|X_v|\leq \kappa$ for all
$v$. Then there is cartesian $R$-submodule $M'\subseteq M$ with
$M'(v)$ pure for all $v$, with $X_v\subseteq M'(v)$ for all $v$
and such that $|M'|\leq \kappa$.
\end{prop}
{\bf Proof. } The proof of \cite[Proposition 3.3]{EE} gives a cartesian
$Q$-submodule $M'$ of $M$ satisfying the desired properties. But then by the
previous comment, as $M$ is cartesian, $M'$ will be also cartesian $R$-module.
$\Box$.

\medskip
\begin{defi}
A cartesian $R$-submodule $M'$ of an $R$-module $M$ is said
to be pure whenever $M'(v)$ is a pure $R(v)$-submodule of $M(v)$,
for every vertex $v\in V$.
\end{defi}

\medskip

\begin{cor}\label{cor1}
There exists an
infinite cardinal $\kappa$ such that every
cartesian $R$-module $M$ is the sum of its quasi-coherent
$R$-submodules of type $\kappa$.

\end{cor} 

{\bf Proof. } Let $M$ be any cartesian $R$-module and take an
element $x\in M$. Then, by Proposition \ref{nocont} we find a
(pure) cartesian $R$-submodule $S_x$ of $M$ with $|S_x|\leq
\kappa$ and $x\in S_x$. Thus $M=\sum_{x\in M} S_x$. $\Box$

\medskip
As a consequence of this we have that $\cart$ is a Grothendieck category
whenever $R$ is a flat presheaf of rings,
for if we take a set $Z$ of representatives of
cartesian modules with cardinality bounded by $\kappa$, it is
immediate that the single cartesian $R$-module $\oplus_{S\in
Z}S$ generates the category of quasi-coherent $R$-modules. 

Now if we focus on particular instances of small categories we have the
following significant consequences. 
\begin{cor}
Let $(X,\O_X)$ be any arbitrary scheme. Each
quasi-coherent sheaf can be written as a continuous chain of pure
quasi-coherent subsheaves of type $\kappa$. Thus $\Qco(X)$ is a
Grothendieck category.

\end{cor}
{\bf Proof.} We let $\C$ consisting of all the affine open $U\subseteq X$. Then
the inclusion between affine open subsets defines a canonical structure of a
partially ordered category on $\C$. Now we let $R$ be the structure sheaf
$\O_X$. Then it is standard that $\cart$ and $\Qco(X)$ are equivalent
categories. So the result will follow from Corollary \ref{cor1}. $\Box$

\medskip\par\noindent
{\bf Remark.}
The previous proof also clarifies a possible misunderstanding on
\cite[Section 2]{EE}. There, the reader may wrongly
think that we are considering the {\it free} category on the affine open
subsets of the scheme $X$ to establish our equivalent category $\mathcal C$.
This is obviously not true, and the gap is easily fixed by saying that we were
assuming the compatibility
condition on our representations there to get the desired equivalence. To be
precise we are just claiming that
$\Qco(X)$ and $\cart$ (or $\C$ in that section) are equivalent.

\medskip\par
Our second application goes back to Artin stacks (cf. \cite{olsson}) and
Deligne-Mumford stacks.

\begin{cor}
Let $\mathcal X$ be a Deligne-Mumford stack. Then the category $\Qco(\mathcal
X)$ is a Grothendieck category. In particular it is locally presentable and has
arbitrary products.
\end{cor}
{\bf Proof.} We take $\C$ as the small subcategory of the iso classes of the
category of affine schemes that are \'etale over $\mathcal X$ (such small
subcategory must exist as the iso classes of such schemes form a set, as
\'etale morphisms are of finite type). Then $\cart$ is equivalent to
$\Qco(\mathcal{X})$.
$\Box$

\begin{cor}
Let $\mathcal X$ be an algebraic stack with a flat sheaf of rings $\mathcal A$.
Then the category $\Qco(\mathcal X)$ is a Grothendieck category. In particular
it is locally presentable and has arbitrary products.
\end{cor}
{\bf Proof.} In this case we consider $\C$ to be the category of affine
schemes smooth
over $\mathcal X$ and $R$ as the sheaf of rings $\mathcal A$. Then $\cart$ is
equivalent to the category of quasi--coherent sheaves on $\mathcal X$. $\Box$

\section{Flat covers and cotorsion envelopes}

\begin{teo}\label{main}

Let $\C$ be any small category and let $R$ be a flat presheaf over $\C$. Then
$\cart$ admits flat covers and cotorsion envelopes.
\end{teo}
{\bf Proof.} Use the same proof as in Theorem 4.1 of \cite{EE}.  $\Box$

\begin{cor}
If $X$ is any scheme, the category $\Qco(X)$ admits flat covers and
cotorsion envelopes.
\end{cor}
{\bf Proof.} Let $M$ be a quasi-coherent sheaf over $X$. We have $\Qco(X)$
equivalent                                                  
to the category $\cart$ of cartesian $R$-modules where $\C$ comes from all
open affine subsets of $X$ and where $R$ comes from the
structure sheaf ${\O}_X$. From the definition of the flat
objects in the two categories we see that the equivalence functors
(in both directions) preserve flatness. Then also since the
functors are clearly additive and exact we get that the property
of being cotorsion is also preserved (since cotorsion is defined
in terms of the splitting of certain short exact sequences). So
then thinking of $M$ as a cartesian $R$-module, $M$ has a
flat cover and a cotorsion envelope in the category of cartesian $R$-modules.
The equivalence then gives the desired
flat cover and cotorsion envelope of $M$ in $\Qco(X)$.$\Box$

\medskip\par
Similarly we have the following consequence:
\begin{cor}
 
Let $\mathcal X$ be a Deligne-Mumford stack (or an Artin stack). The
category $\Qco(\mathcal X)$ admits flat covers and
cotorsion envelopes.
\end{cor}

\end{document}